\documentstyle[12pt]{article}

  \newcommand{\const}{\rm const}
  \newcommand{\Var}{\rm Var}
  
  \newcommand{\Law}{\rm Law}

  \newcommand{\vraisup}{\rm vraisup}
  \newcommand{\Dom}{\rm  Dom}

  \newcommand{\argmin}{\rm argmin}
  
  \newcommand{\Cov}{\rm  Cov}
 \newcommand{\Sub}{\rm  Sub}
 \newcommand{\StSub}{\rm  StSub}
 \newcommand{\mod}{\rm mod}

\vspace{3mm}

\begin{document}

   \begin{center}

 {\bf Adaptive  smoothness of function estimation }\\

\vspace{4mm}

{\bf in the three classical problems of the } \\

\vspace{4mm}

{\bf non-parametrical statistic} \\

\vspace{4mm}

 {\bf M.R.Formica,  E.Ostrovsky and L.Sirota.}

 \end{center}

\vspace{4mm}

 Universit\`{a} degli Studi di Napoli Parthenope, via Generale Parisi 13, Palazzo Pacanowsky, 80132,
Napoli, Italy. \\

e-mail: mara.formica@uniparthenope.it \\

\vspace{3mm}

Department of Mathematics and Statistics, Bar-Ilan University, \\
59200, Ramat Gan, Israel. \\

e-mail: eugostrovsky@list.ru\\

\vspace{3mm}

\ Department of Mathematics and Statistics, Bar-Ilan University, \\
59200, Ramat Gan, Israel. \\

e-mail: sirota3@bezeqint.net \\

\vspace{4mm}

\begin{center}

  {\bf Abstract} \par

\end{center}

\vspace{3mm}

\hspace{3mm} We offer in this short report the so-called adaptive  functional smoothness estimation
in the Hilbert space norm sense in the three classical problems of non - parametrical statistic: \par

\begin{center}

 \ {\it   regression, density and spectral (density) function measurement (estimation).} \par

\end{center}

\vspace{4mm}

 \begin{center}

 \ {\sc Key words and phrases:}

 \vspace{3mm}

 \end{center}

 \hspace{3mm} Three classical problem of non - parametrical statistic, adaptive estimation, probability, ordinary
and spectral  density, regression and spectral function, covariation, trigonometrical orthonormal sequence,
 centered ergodic stationary Gaussian sequence, expectation, variance, index of smoothness,
generating and moment generating functions,  empirical and ordinary Fourier coefficients, smoothness of function, expectation and variance,
 ordinary and Young - Fenchel (Legendre) transform, Central Limit Theorem (CLT), uniform norm and estimate, Kramer's condition,
 estimation (measurement), class of smoothness,  convergence, Fourier series and coefficients, Lebesgue - Riesz and
 Grand Lebesgue Spaces (GLS) and norms, distribution and tail of distribution,   confidence interval. \\

\vspace{4mm}

\section{Introductions. Notations. Statement of problem.}

\vspace{4mm}

 \hspace{3mm} Let us consider the following three classical problem of the non - parametric statistics. \par
{\it  Our object of estimation  is a numerical valued (measurable) function} $ \ f = f(x), \ x \in [0,1]. \ $  \\

 \vspace{3mm}

 \ {\bf A. Estimation of the regression function.} \par

\vspace{3mm} Model: given the variables $ \ \{y_i\}, \ i = 1,2,\ldots,n; \ $  such that

\vspace{3mm}

\begin{equation} \label{Regr}
y_i = f(x_i) + \xi_i, \ i = 1,2,\ldots, n;  \ x_i = x_i(n)= i/n; \ n = 2,3,\ldots;
\end{equation}

\vspace{3mm}

where $ \  \{\xi_i \} \ $ are centered independent identically distributed  (c, i,i.d.) random variables (errors of measurement).  \par
 Hereafter the conditions on the  $\ f, \ \{\xi_i \} \ $  will be confirmed. \par

\vspace{3mm}

 \ {\bf B. Density estimation.   } \par

\vspace{3mm}

 \ Datum: a sample $ \  \{\xi_i\}, \ i = 1,2,\ldots,n  \ $ as ordinary a sequence of  \ i, i. d \ random variables  having a common density $ \ f = f(x): \ $

\vspace{3mm}

$$
{\bf P} (\xi_i \in (a,b)) = \int_a^b f(x) \ dx; \ 0 \le a < b \le 1.
$$

\vspace{3mm}

 \ {\bf C.  Spectral density (function) measurement.} \par

 \vspace{3mm}

 \ Datum: a Gaussian distributed random centered stationary sequence $ \ \{ \xi_i \}, \ i = 1,2,\ldots,n \ $ with spectral density $ \ f = f(x), \ x \in [0,1]: \ $

 \vspace{3mm}

 $$
\Cov(\xi_i,\xi_j) = {\bf E} \xi_i \xi_j  = \sqrt{2} \ \int_0^1 \cos(2 \pi (i - j) x) \ f(x) \ dx.
 $$

 \vspace{5mm}

 \ Introduce the ordinary orthonormal  trigonometrical sequence of functions $ \ \{\phi_k(x) \},  x \in [0,1] \ $  as follows: \ $ \ \phi_1(x) = 1; \ $ \par

$$
\phi_{2l + 1}(x) \stackrel{def}{=}\sqrt{2} \sin(2 \pi l x); \hspace{3mm}  \phi_{2l}(x) \stackrel{def}{=}\sqrt{2} \cos(2 \pi l x), \ l = 1,2,\ldots.
$$

\vspace{3mm}

\hspace{3mm} {\sc We suppose hereafter in particular that the unknown (estimated) function belongs to the usually Hilbert space } $ \ L_2[0,1]: \ $

\vspace{3mm}

$$
|| \ f \ ||^2 L_2[0,1] \stackrel{def}{=}  \int_0^1 f^2(x) dx < \infty,
$$
where integral is understood in the Riemann sense.\par

\vspace{3mm}

 \ Therefore it may be represented by the classical Fourier series convergent in $ \ L_2[0,1]  \ $ sense

$$
f(x) = \sum_{k=1}^{\infty} c_k \ \phi_k(x),
$$
where  $ \  \sum_{k=1}^{\infty} c_k^2 < \infty,  \ $ and hence  $ \ \lim_{N \to \infty} \rho(N) = 0, \ $ where

\vspace{3mm}

\begin{equation} \label{rho def}
\rho(N) = \rho[f](N) \stackrel{def}{=} \sum_{k = N + 1}^{\infty} c_k^2.
\end{equation}

\vspace{3mm}

 \ Obviously,

$$
  \rho(0) = \sum_{k = 1}^{\infty} c_k^2 = \int_0^1 f^2(x) dx = ||f||^2 L_2[0.1].
$$

\vspace{4mm}

  \hspace{3mm} {\sc  The class of smoothness in the Hilbert space sense norm  $ \ L_2[0,1] \otimes L_2[\Omega] \ $  for the estimated
function is by definition closely related  with speed  of a convergence } $ \ \rho[f](N) \to 0 \ $  {\sc as}   $ \ N \to \infty. \ $ \par

\vspace{3mm}

 \hspace{3mm} For example,

$$
||f||C[0,1] = \max_{x \in [0,1]} |f(x)| \le \sum_{l=0}^{\infty} 2^{l + 1/2} \ \rho^{1/2}[f] \left(2^l \right),
$$
if the series here converges; in addition, the function $\ f(\cdot) \ $  is continuous (Nikol'skii inequality). \par

\vspace{4mm}

 \ Recall that correspondent  (asymptotical) unbiased consistent  as $ \ n \to \infty \ $ estimates of the Fourier coefficients $ \ \hat{c}_k  = \hat{c}_k(n) \ $
 has the following forms.  \par

 \vspace{3mm}

 \ {\bf A. For the estimation of the regression function}

\vspace{3mm}

\begin{equation} \label{regr coef}
\hat{c}_k := n^{-1} \sum_{i=1}^n y_i \phi_k(x_i(n)).
\end{equation}

\vspace{3mm}

 \ {\bf B. For the problem  of density estimation}

\vspace{3mm}

\begin{equation} \label{regr coef}
\hat{c}_k := n^{-1} \sum_{i=1}^n  \phi_k(\xi_i).
\end{equation}

\vspace{3mm}

 \ {\bf C. For the spectral density estimation}

\vspace{3mm}

\begin{equation} \label{regr coef}
\hat{c}_k := n^{-1} \sum_{i=1}^{n - k}  \xi_i \xi_{i + k}.
\end{equation}

\vspace{4mm}

\section{Preliminary estimate.}

\vspace{3mm}

 \ Define the following important variable

 \vspace{3mm}

\begin{equation} \label{Tchentzsov}
A(n,N) \stackrel{def}{=} \frac{N}{n} +\rho(N), \ N = 1,2,\ldots,n.
\end{equation}

\vspace{3mm}

 \hspace{3mm} N.N.Tchentzov in \cite{Tchentzov1}, \cite{Tchentzov2} offered the following so - called  "projection"
 estimates for the function $ \ f \ $ in all the considered three problems.

\vspace{3mm}

\begin{equation} \label{Err1}
\hat{f}_{n,N}(x) \stackrel{def}{=} \sum_{k=1}^N \hat{c}_k(n)\phi_k(x).
\end{equation}

 \vspace{3mm}

 \ It is no hard to calculate

\vspace{3mm}

\begin{equation} \label{non adapt}
{\bf E} ||\hat{f}_{n,N} - f||^2 L_2[0,1] = A(n,N).
\end{equation}

\vspace{3mm}

\ Therefore the  optimal choice of the numbers of harmonics $ \ N^* = N^*(n) \ $ based on the  offered method
can be implemented by the formula

\begin{equation} \label{non - ad}
N^* = N^*(n) = \argmin_{N \in [1,n]} A(n,N);
\end{equation}

\vspace{3mm}

 and correspondingly $ \ A^* (n) \stackrel{def}{=}  \min_{N \in [1,n]} A(n,N) = A(n, N^*(n));    \ $

\vspace{3mm}

\begin{equation} \label{non ad optimal}
\min_{N \in [1,n]} {\bf E}  \left[ \ ||\hat{f}_{n,N} - f ||L_2[0,1] \ \right]^2  =  {\bf E}  \left[ \ ||\hat{f}_{n,N^*(n)} - f ||L_2[0,1] \ \right]^2 = A^*(n).
\end{equation}

\vspace{3mm}
\ N.N.Tchentzov in \cite{Tchentzov1}, \cite{Tchentzov2} proved in particular that the estimate $ \ \hat{f}_{n,N^*(n)} \ $  is asymptotically optimal
 as $ \ n \to \infty \ $ on the wide class of Banach functional spaces containing the unknown function $ \ f(\cdot). \ $ \par
 \ {\it But this  estimates  are non adaptive,}  as long as the behavior of the function $ \ N \to \rho[f](N) \ $
 as $ \ N \to \infty \ $ is in general case unknown. The {\it adaptive} modification of these estimations is offered and investigated e.g. in
\cite{Bobrov}, \cite{Bobrov2},  \cite{Ostrovsky1}, chapter 5, section 5.13. \par
 \ In detail, introduce the following important variables (statistics)

 \vspace{3mm}

\begin{equation} \label{tau}
\tau_N = \tau_N(n)  \stackrel{def}{=} \sum_{j=N+1}^{2N} \hat{c}_j^2; \hspace{3mm}  \tilde{N} = \tilde{N}(n) = \tilde{N}[\vec{\xi}](n) \stackrel{def}{=} \argmin_{N \in [1,n]}\tau_N(n),
\end{equation}

\vspace{3mm}

\begin{equation}  \label{min tau}
\tau^*(n) \stackrel{def}{=} \min_{N \in [1,n]} \tau_N(n);
\end{equation}

\vspace{3mm}

so that

\vspace{3mm}

\begin{equation} \label{min tau}
\tau^*(n) = \tau_{\tilde{N}(n)}(n).
\end{equation}

\vspace{4mm}

 \hspace{3mm}  {\it We suppose hereafter that
 the estimated function  in all the considered problems}
 $ \ f(\cdot) \ $ {\it satisfies the following condition}

 \vspace{3mm}

\begin{equation} \label{gamma condit}
\exists \gamma = \const \in (0,1), \ \Rightarrow \sup_{N \ge 2} \rho[f](2N)/\rho[f](N) \le \gamma,
\end{equation}

 \vspace{3mm}

 or more generally

 \vspace{3mm}

\begin{equation} \label{cond gamma}
\overline{\lim}_{N \to \infty} \rho[f](2N)/\rho[f](N) < 1,
\end{equation}

 \vspace{3mm}

then the new estimate (already adaptive!) has a form

\vspace{3mm}

\begin{equation} \label{ad estim}
\tilde{f}_n(x) =  \tilde{f}[\vec{\xi}]_n(x)   \stackrel{def}{=} \sum_{k=1}^{ \tilde{N}(n)} \hat{c}_k(n) \phi_k(x)
\end{equation}

\vspace{3mm}

is asymptotical optimal as $ \ n \to \infty. \ $  In detail, it is proved in \cite{Bobrov},  \cite{Ostrovsky1},
chapter 5, section 5.13  that as $ \ n \to \infty \ $

\vspace{3mm}

\begin{equation} \label{adap  optimal}
   {\bf E}  \left[ \ ||\tilde{f}_n - f ||L_2[0,1] \ \right]^2 \asymp A^*(n).
\end{equation}

\vspace{3mm}

 \ Moreover, \ $ \ \lim_{n \to \infty} \tilde{N}(n)/N^*(n) = 1; \ $ \\

\vspace{3mm}

\begin{equation} \label{moreover}
{\bf E} \tau_N^* (n) \asymp A(n,N); \hspace{3mm} \lim_{n \to \infty} \frac{\tau_N^*(n)}{{\bf E} \tau_N(n)} =  1
\end{equation}

\vspace{3mm}

in probability.  Therefore, there is some reason to believe that the following {\it adaptive} estimate is consistent, of course
under appropriate  natural conditions

\vspace{3mm}

\begin{equation} \label{f ad estim}
\tilde{f}_n(x) \stackrel{def}{=} \sum_{k=1}^{\tilde{N}} \hat{c}_k \ \phi_k(x) \stackrel{P}{\to} f(x).
\end{equation}

\vspace{3mm}

in the $ \ L_2[0,1] \ $ sense;  and moreover as $ \ n \to \infty \ $

\vspace{3mm}

\begin{equation} \label{key relation}
\tau^*(n) \asymp  \ \min_{N \in [1,n]} \left( \ \frac{N}{n} + \rho(N) \ \right).
\end{equation}

\vspace{3mm}
 \ The relation   (\ref{f ad estim}) is proved in particular in \cite{Bobrov},  \cite{Ostrovsky1}, chapter 5, section 5.13.
 The relation   (\ref{key relation}), which will be proved further, may be used, in particular for the estimation the value $ \ \rho(N), \ $
 as well.\par

\vspace{3mm}

 \ Since the case when $ \ \exists N_0 \in (2,3, \ldots) \ \Rightarrow \rho(N_0) = 0, \ $ i.e. when the estimated function is trigonometrical
 polynomial, is trivial, we suppose further that $ \ \forall N \ \Rightarrow \ \rho(N) > 0. \ $ \par
  \ The variable $ \ \rho(N) \ $ may be named as a  {\it  index of smoothness} of the estimated function. \par

\vspace{4mm}

\hspace{3mm} {\bf  Our aim in this short report is to offer the consistent } as $ \ n \to \infty \ $ {\bf adaptive estimates for the class of smoothness
for considered function } $ \ \rho(N) \ $  {\bf in all the three statements of problem. } \par
 \ We establish also the construction of a confidence  region for this function, asymptotical or not. \par

\vspace{3mm}

  \hspace{3mm} {\it Preliminary estimates.} \par

 \vspace{3mm}

 \ Note that it is reasonable   to suppose that as $ \ n \to \infty \ $

\vspace{3mm}

 \begin{equation} \label{reasonable}
\tau(n,N) \asymp A(n,N) = \frac{N}{n} + \rho(N), \ N \in [1,n]; \ N \to \infty.
 \end{equation}

\vspace{3mm}

 \ Therefore, one can offer and investigate the following {\it consistent} as $ \ n \to \infty \ $ estimate  $ \ \hat{\rho}_n(N) \ $ for the
 {\it smoothness function} $ \ \rho(N): \ $

\vspace{3mm}

 \hspace{3mm} {\bf Definition 2.1.}

 \begin{equation} \label{prelim estim}
 \hat{\rho}(N) =  \hat{\rho}_n(N) \stackrel{def}{=} \tau(n,N) - \frac{N}{n}.
 \end{equation}

\vspace{4mm}

 \ {\bf Remark 2.1.} Note that the optimal adaptive estimations of the function $ \ f \ $ in the uniform
 norm was investigated in works \cite{Golubev Lepsky Levit}, \cite{Bobrov},\cite{Bobrov2}.\par

\vspace{3mm}

\section{Asymptotical analysis.}

\vspace{3mm}

\hspace{3mm} Here $ \ n,N \to \infty \ $ and hence we apply the classical Central Limit Theorem (CLT).
\  In all these statements holds true the following expression for the empirical  Fourier coefficients

 \vspace{3mm}

\begin{equation} \label{key expression}
\hat{c}_k(n) \sim c_k +  \sigma \cdot \frac{\delta_k}{\sqrt{n}}, \hspace{3mm}  k,n \to \infty,
\end{equation}

 \vspace{3mm}

where $ \ \{\delta_k\} \ $ are centered asymptotically independent
 non - zero  normal identical distributed random variables: \ $ \ \Law \{\delta_k \} = N(0,1), \ 0 < \sigma = \const < \infty. \ $
 For instance,  it is easily to calculate for the first discussed problem  \ {\bf A}

$$
\sigma^2 = \Var (\xi_1),
$$
  if of course it is finite. \ In the  problem {\bf B}

 $$
   \sigma = 1, \ \Var (\delta_k) \sim 1, \  k \to \infty.
 $$

\vspace{3mm}

 \ Finally, in the third problem {\bf C}

$$
\sigma^2 = 2 c^2_0  + \sum_{j=1}^{\infty} c^2_k.
$$

\vspace{4mm}

 \ One can assume further  for calculations  without loss of generality $ \ \sigma = 1. \ $ \par
 \ We have

$$
\tau_n(N) = \sum_{k=N+1}^{2N}\left[\hat{c}_k \right]^2 = \sum_{k=N+1}^{2N} c_k^2  + \frac{2}{\sqrt{n}}  \sum_{k=N+1}^{2N} c_k \delta_k + \frac{1}{n}\sum_{k=N+1}^{2N} \delta_k^2 =
$$

$$
A(N,n) +  \frac{2}{\sqrt{n}}  \sum_{k=N+1}^{2N} c_k \delta_k + \frac{1}{n}\sum_{k=N+1}^{2N} \left(\delta_k^2 - 1 \right) = A(N,n)  + \Sigma_1 + \Sigma_2;
$$

where

$$
\Sigma_1 =  \frac{2}{\sqrt{n}}  \sum_{k=N+1}^{2N} c_k \delta_k,
$$

$$
\Sigma_2 = \frac{1}{n}\sum_{k=N+1}^{2N} \left(\delta_k^2 - \sigma^2 \right).
$$

 \ We  have obviously

$$
||\Sigma_1||L_2(\Omega) \le 2 n^{-1/2} \ \sqrt{\rho(N)}; \ ||\Sigma_2||L_2(\Omega) \le C  \ n^{-1} \ \sqrt{N};
$$

$$
 2 n^{-1/2} \ \sqrt{\rho(N)} \le \frac{1}{n} + \rho(N);
$$

so that as $ \ n \to \infty  \ \Rightarrow \ $

$$
\lim  \frac{\tau_n(N)}{A(N,n)} =1,  \  N \in [1,n]
$$
in probability and in the Hilbert norm  $ \   L_2(\Omega, {\bf P}) \ $ sense. \par

\vspace{3mm}

  \  Therefore $ \tau_n(N) \sim A(N,n);  \   \tau_n(N) - \frac{N}{n} \sim A(N,n) -  \frac{N}{n}.\   $  To summarize: \par

 \vspace{3mm}

 \ {\bf Proposition 3.1.} We conclude under formulated  notations and conditions as $ \ n \to \infty \ $ for all the values $ \ N \ $
 from the integer set $ \ [1, n - 1 ]. \ $\\

\vspace{3mm}

\begin{equation} \label{asymp est}
  \lim \frac{\hat{\rho}_n(N)}{\rho(N)} = 1 \ (\mod \ {\bf P}).
\end{equation}

 \vspace{3mm}

 \ Therefore, the variable $ \ \hat{\rho}_n(N) \ $ is consistent as $ \ n \to \infty \ $ estimate for the one $ \ \rho_n(N) \ $
for arbitrary {\it fixed} value $ \ N; \ N \in [2, n-1]. \ $ But appears an open question about {\it asymptotical} confidence interval.\par

\vspace{3mm}

 \ We have the following expression

\vspace{3mm}

 \begin{equation} \label{expression}
 \hat{\rho}_n(N) = \rho_n(N) + \frac{2}{\sqrt{n}} \sum_{N+1}^{2N} c_k \delta_k + \frac{1}{n} \sum_{k= N+1}^{2N} (\delta_k^2 - 1) =: \rho_n(N) + S_1 + S_2.
 \end{equation}

\vspace{3mm}

 \ Denote

$$
\Theta(n,N,\rho(N)) := \sqrt{ \ \frac{4}{n} \rho(N) + \frac{9}{n^2} N \ \ };
$$

then we deduce by virtue of  the classical CLT the random variable  $ \ \hat{\rho}_n(N) \ $ has approximately Gaussian distribution: as $ \ n \to \infty \ $

$$
\Law \left[ \ \hat{\rho}_n(N) \ \right] \sim  N(\rho_n(N), \ \Theta^2(n,N, \rho(N))).
$$
 \ This fact allow us to build the confidence interval for the value  $ \ \rho_n(N) \ $   as follows. Let  the value $ \ \alpha \in (0,1/2) \ $
 be given. Define a positive number $ \ v = v(\alpha) \ $ as usually

 $$
 (2 \pi)^{-1/2} \int_{-v}^{v} \exp(-x^2/2) dx = \alpha;
 $$

then with probability $ \ \approx \alpha \ $ there holds

\vspace{3mm}

\begin{equation} \label{Appr conf interval}
|\hat{\rho}_n(N) - \rho_n(N)|/\Theta(n,N,\rho(N)) \le v(\alpha).
\end{equation}

\vspace{4mm}

\section{Auxiliary buildings.}

\vspace{4mm}

\begin{center}

 {\sc Grand Lebesgue Spaces.}

\vspace{4mm}

\end{center}

 \hspace{3mm}  We recall here for reader convenience some known definitions and  facts  from the theory of Grand Lebesgue Spaces (GLS).
    \ Let $ \ \psi = \psi(p), \ p \in [1,b) \ $    where $ \ b = \const \in (1,\infty] \ $ be positive measurable numerical valued
 function, not necessary to be finite in the boundary  point $ \ p = b-0 \, $ such that $ \ \inf_{p \in [1,b)} \psi(p) > 0. \ $ The set of
 all such a functions will be denoted as $ \ G \Psi(b); \ $  define also for certain numerical values $ \ b > 1 \ G\psi(b) \ $
some function from the set  $ \ G\Psi(b): \  G\psi(b) \in G\Psi(b), \ $ and put also

$$
G\Psi := \cup_{1 < b \le \infty} G \Psi(b).
$$

\vspace{3mm}

 \  For instance

  $$
    \psi_m(p) := p^{1/m}, \ m = \const > 0, \ p \in [1,\infty)
  $$
   or

$$
   \psi^{b; \beta}(p) :=  (b-p)^{-\beta}, \ p \in [1,b), \ 1 < b \le \infty, \ \beta = \const > 0.
$$

 \ The case $ \ m = 2 \ $  correspondent to the classical {\it subgaussian} space: $ \ ||\zeta||G\psi_2 = ||\zeta||\Sub. \ $

\vspace{3mm}

 \ The set of all such a {\it generating} functions will be denoted  by $ \   \Psi = G\Psi(b) = \{\psi(\cdot)\}, \ $ as well as a
 correspondent Banach space. \par

 \vspace{3mm}

  \hspace{3mm} {\bf Definition 4.1.}  The (Banach) Grand Lebesgue Space (GLS)    $  \ G \psi  = G\psi(b),  $
    consists on all the real (or complex) numerical valued measurable functions (random variables!)
   $   \  f: \Omega \to R \ $  defined on whole our  probability space $ \  \Omega \ $ and having a finite norm

\vspace{3mm}

 \begin{equation} \label{norm psi}
    || \ f \ || = ||f||G\psi \stackrel{def}{=} \sup_{p \in [1,b)} \left[ \frac{||f||_p}{\psi(p)} \right].
 \end{equation}

 \vspace{4mm}

  \ As usually,

$$
||f||_p \stackrel{def}{=} \left[ \ {\bf E} |f|^p  \ \right]^{1/p}, \ 1 \le p < \infty.
$$

\vspace{3mm}

 \ The function $ \  \psi = \psi(p) \  $ is named as  the {\it  generating function } for this space. \par

  \ If for instance

$$
  \psi(p) = \psi^{(r)}(p) = 1, \ p = r;  \  \psi^{(r)}(p) = +\infty,   \ p \ne r,
$$
 where $ \ r = \const \in [1,\infty),  \ C/\infty := 0, \ C \in R, \ $ (an extremal case), then the correspondent
 $ \  G\psi^{(r)}(p)  \  $ space coincides  with the classical Lebesgue - Riesz space $ \ L_r = L_r(\Omega). \ $ \par

\vspace{4mm}

 \ These spaces are investigated in many works, e.g. in \cite{Ahmed Fiorenza Formica at all}, \cite{Ermakov etc. 1986},
 \cite{Fiorenza1},   \cite{Fiorenza-Formica-Gogatishvili-DEA2018}, \cite{Fiorenza4},    \cite{Iwaniec2}, \cite{KozOs},
\cite{liflyandostrovskysirotaturkish2010}, \cite{Ostrovsky1}  - \cite{Ostrovsky4}, \cite{Yudovich1}, \cite{Yudovich2} etc. They are applied for example in the
theory of Partial Differential Equations \cite{Fiorenza-Formica-Gogatishvili-DEA2018},
 \cite{Fiorenza4}, in the theory of Probability  \cite{Ermakov etc. 1986}, \cite{Ostrovsky3}  - \cite{Ostrovsky4}, in Statistics \cite{Ostrovsky1}, chapter 5,
theory of random fields   \cite{KozOs}, \cite{Ostrovsky4}, in the Functional Analysis \cite{Ostrovsky1}, \cite{Ostrovsky2}, \cite{Ostrovsky4} and so one. \par

 \  These spaces are rearrangement invariant (r.i.) Banach functional  (complete) spaces. They do
  not coincides  in general case with the classical Banach rearrangement functional  spaces: Orlicz, Lorentz, Marcinkiewicz  etc., as well; see
  \cite{liflyandostrovskysirotaturkish2010}, \cite{Ostrovsky2}.\par

   \ The belonging of certain  measurable function (random variable) $ \ f: \Omega \to R \ $ to some $ \ G\psi \ $ space   is closely related with  its {\it tail} function

 $$
 T_f(t) \stackrel{def}{=} {\bf P}(|f| \ge t), \ t \ge 0,
 $$
 behavior   when $ \ t \to \infty,  \ $  see  \cite{Ahmed Fiorenza Formica at all}, \cite{KozOs}, and so one. \par

 \hspace{3mm} In detail, let the non - zero random variable $ \ \upsilon \ $ belongs to the certain Grand Lebesgue Space $ \ G\psi \ $ and, let for definiteness
 $ \ ||\upsilon||G\psi = 1. \ $ Define an auxiliary  functions

$$
h(p) = h[\psi](p) := p \ln \psi(p), \ h^*[\psi](v) := \sup_{p \in \Dom[\psi]}(p v - h(p)) \ -
$$
 the so - called Young - Fenchel, or Legendre transform for  the function $ \  h(p). \ $ The symbol $ \ \Dom[\psi] \ $
 denotes the domain of definition (and finiteness) for the function $ \ \psi(\cdot). \ $ \par

\vspace{3mm}

  \  Holds true

\begin{equation} \label{Young Fen}
T_{\upsilon}(t) \le \exp(-h^*(\ln t)), \ t \ge e.
\end{equation}
 \ Note that the inverse conclusion is true, of course under suitable  natural conditions, see  \cite{KozOs}, \cite{Ostrovsky1}. \par

\vspace{3mm}

 \ Further let $ \ \xi \in G\psi = G\psi(b),  \ ||\xi||G\psi  =:\Delta < \infty, \ \psi \in G\Psi, \ m = \const, b \ge m > 1. \ $ There holds

$$
{\bf E} |\xi|^p \le \Delta^p \ \psi^p(p), \ p \in[1,b).
$$

 \ Define a new generating function as follows

$$
\psi^{(m)}(p) := \psi^m(mp), \ 1 \le p \le b/m.
$$

  \ We conclude

$$
{\bf E} |\xi|^{mp} \le \psi^{mp}(mp) \cdot \Delta^{mp},
$$
or equally

\vspace{3mm}

\begin{equation} \label{degree m}
|| \ |\xi|^m \ || G\psi^{(m)} \le (||\xi|| G\psi)^m.
\end{equation}

\vspace{3mm}

 \ For instance

\vspace{3mm}

\begin{equation} \label{degree 2}
|| \ \xi^2 \ || G\psi^{(2)} \le (||\xi|| G\psi)^2.
\end{equation}

\vspace{5mm}

 \begin{center}

 \ {\sc B spaces.}

 \end{center}

\vspace{3mm}

 \ Let now  $ \lambda_0 = \const \in (0,\infty]$ and let $ \phi = \phi(\lambda) $ be an
even strong convex function in the set  $ \ (-\lambda_0, \lambda_0), \ $ which takes only
positive values for non - zero argument, twice continuously differentiable; briefly
$ \ \phi = \phi(\lambda) \ $ is a Young-Orlicz function, and such that

\vspace{3mm}

\begin{equation}\label{Young-Orlicz function}
\phi(0) = 0, \ \ \phi'(0) = 0, \ \  \phi^{''}(0) \in (0,\infty).
\end{equation}

\vspace{3mm}

 \   We denote the set of all these Young - Orlicz function as  $ \ \Phi =\Phi[\lambda_0]: \ \Phi = \{ \phi(\cdot) \ \}. $  \par

\vspace{4mm}

 \ {\bf Definition  4.2.} \par

\vspace{3mm}

 \ Let $\phi\in \Phi$. We say that the centered random variable $\xi$
belongs to the space $B(\phi)$  if there exists a constant $\tau
\geq 0$ such that
\begin{equation}\label{spaceB}
\forall \lambda \in (-\lambda_0, \lambda_0) \ \Rightarrow {\bf E}
\exp(\pm \lambda \ \xi) \le \exp(\phi(\lambda \ \tau)).
\end{equation}
 The minimal non-negative value $\tau$ satisfying (\ref{spaceB}) for
arbitrary  $\lambda \in (-\lambda_0, \ \lambda_0)$ is named $B(\phi)$-norm
of the variable $\xi$  and we will write

\vspace{3mm}

\begin{equation}\label{Bnorm}
||\xi||B(\phi) \stackrel{def}{=}\inf \{\tau, \ \tau \ge 0 \ : \ \forall
\lambda \in (-\lambda_0, \lambda_0) \ \Rightarrow {\bf E} \exp(\pm
\lambda \ \xi) \le \exp(\phi(\lambda \ \tau)) \} .
\end{equation}

\vspace{3mm}

 \ For instance if $\phi(\lambda)=\phi_2(\lambda) := 0.5 \ \lambda^2, \
\lambda \in \mathbf{R}$, the r.v. $\xi$ is \emph{subgaussian} and in
this case we denote the space $B(\phi_2)$ with $\Sub$. Namely we
write $\xi \in \Sub$ and

$$
||\xi||\Sub \stackrel{def}{=} ||\xi||B(\phi_2).
$$

 \vspace{3mm}

  \ Obviously, arbitrary centered (non - zero) Gaussian distributed r.v. $ \ \zeta \ $  is subgaussian. Moreover, it
is {\it strictly subgaussian}; this means by definition that $ \ ||\zeta||\Sub = \{ \Var [\zeta] \}^{1/2}. \ $ Another example of
such a variables are so - called {\it Rademacher's variables:}

$$
{\bf P} (\zeta = 1) = {\bf P}(\zeta = - 1) = 1/2.
$$

 \ Here $ \ {\bf E}  \exp(\lambda \ \zeta) = \cosh \lambda \le \exp(0.5 \lambda^2), \ \lambda \in R. \  $ \par

\vspace{3mm}

 \ {\it Tail behavior for these r.v.} \ Let the non - zero r.v. $ \ \xi \ $ belongs to some  $ \ B\phi \ $ space:
$ \ ||\xi||B\phi = v \in (0,\infty). \ $ Then

\vspace{3mm}

\begin{equation} \label{tail Bphi}
{\bf P} ( \ |\xi| > t \ ) \le 2 \exp( \ - \phi^*(t/v) \ ), \ t > 0,
\end{equation}

\vspace{3mm}

 the famous Tchernov's estimate. The inverse conclusion also holds true up to multiplicative constant, of course under natural
 appropriate conditions.\par

\vspace{3mm}

 \hspace{3mm} Notice again that if the {\it centered} random variable $ \  \zeta  \ $ is a.s. bounded, then

\vspace{3mm}

\begin{equation} \label{boundedness}
|| \zeta \ ||\Sub \le \vraisup_{\omega \in \Omega} |\zeta(\omega)|.
\end{equation}

\vspace{3mm}

 \  The following example will be used further. Define a new generating function  for these $ \ B \ $ spaces:

\vspace{3mm}

\begin{equation} \label{ups fun}
\upsilon(\lambda) \stackrel{def}{=}  -0.5 \ln(1 - 2 |\lambda|) - |\lambda|; \hspace{3mm}  |\lambda| < 1/2 = \lambda_0.
\end{equation}

\vspace{3mm}

 \ Let as before $ \ \Law(\xi) =  N(0,1);  \  $ and set $ \ \beta := \xi^2 - 1; \ $ then $ \ {\bf E} \beta = 0   \ $ and

$$
{\bf E} \exp(\lambda \ \beta) = \exp  \upsilon(\lambda), \  \lambda \in (-0.5; +0.5);
$$
and of course $ \ {\bf E} \exp(\lambda \ \beta) = +\infty, \ |\lambda| \ge 1/2. \ $  \par
 \ On the other words, $ \ ||\beta||B\upsilon = 1. \ $ \par
 \ More generally,   let $ \ \Law(\xi) = N(0,\Delta^2), \ 0 < \Delta < \infty;  \ $  and put
 $ \ \tau = \xi^2 - \Delta^2. \ $ Then

$$
{\bf E} \exp \left[ \ \lambda \ \tau   \ \right] = \exp \left[\upsilon (\lambda \ \Delta^2) \ \right], \  |\lambda| < 1/(2 \Delta^2),
$$
or equally

$$
 ||\tau||B\upsilon = \Delta^2.
$$

\vspace{4mm}

 \hspace{3mm}  It is known, see  \cite{KozOs}, \cite{Buldygin-Mushtary-Ostrovsky-Pushalsky} that if the random variables $ \ \xi_i \ $ are
independent and subgaussian, then

$$
||\sum_{i=1}^n \xi_i||_{\Sub} \le \sqrt{\sum_{i=1}^n ||\xi_i||^2_{\Sub}}.
$$

 \ At the same inequality holds true in the more general case in the $ \ B(\phi) \ $ norm, when the function
 $ \ \lambda \to \phi(\sqrt{|\lambda|}), \ |\lambda| < \lambda_0  \ $ is convex, see \cite{KozOs}. \par
 \ As a slight corollary: in this case and if in addition the r.v. - s $ \ \{\xi_i \} \ $ are i., i.d., centered,  then

$$
 \sup_{n = 1,2,\ldots}|| n^{-1/2} \sum_{i=1}^n \xi_i||B(\phi) =  ||\xi_1||B(\phi).
$$

\vspace{3mm}

 {\bf Definition 4.3.} \

\vspace{3mm}

 \ The centered r.v. $\xi$ with finite non-zero variance $\sigma^2 :=
\Var (\xi) \in (0,\infty)$ is said to be {\it strictly subgaussian,} and
write $ \ \xi \in \StSub \ $, iff
$$
{\bf E} \exp(\pm \ \lambda \ \xi) \le \exp(0.5 \ \sigma^2 \
\lambda^2), \ \lambda \in R.
$$

\vspace{4mm}

 \ For instance, every centered non-zero Gaussian r.v. belongs to the
space $\StSub.$ The Rademacher's r.v. $\xi$, that is such that
${\bf P}(\xi = 1) = {\bf P}(\xi = - 1) = 1/2$, is also strictly
subgaussian. Many other strictly subgaussian r.v. are represented in
\cite{Buldygin-Mushtary-Ostrovsky-Pushalsky}, \cite{KozOs},  \cite{Ostrovsky1}, chapter 1. \par

\vspace{3mm}

\ It is proved in particular that $ \ B(\phi), \ \phi  \in \Phi, \ $ equipped with the norm
   (\ref{Bnorm}) and under the ordinary algebraic operations, are
Banach rearrangement invariant  functional spaces, which are
equivalent the so-called Grand Lebesgue spaces as well as to {\it Orlicz
exponential spaces.} These spaces are very convenient for the
investigation of the r.v. having an exponential decreasing tail of
distribution; for instance, for investigation of the limit theorem,
the exponential bounds of distribution for sums of random variables,
non-asymptotical properties, problem of continuous and weak
compactness of random fields, study of Central Limit Theorem in the
Banach space, etc. \par

\vspace{3mm}

 \ Let $ \ \phi(\cdot) \in \Phi. \ $  Define the following such a function

\vspace{3mm}

\begin{equation} \label{aux der fun}
\overline{\phi}(\lambda) \stackrel{def}{=} \sup  \left\{ \ \sum_{j=1}^{\infty} \phi(\gamma_j \ \lambda) \ \right\},
\end{equation}

\vspace{3mm}

 \ where $ \ \sup \ $ is calculated over all the numerical sequences $\ \{\gamma_j \} \ $ such that $ \ \sum_j \gamma_j^2 = 1. \ $ \par

\vspace{3mm}

 \ This definition obeys a following sense: if the independent r. v.  $ \   \{\eta_j\} \ $ are normed bounded in the space
  $ \ B(\phi): \ \sup_j ||\eta_j||B \phi  = \theta \in (0,\infty), \ $ then  also when $ \ \sum_j \gamma_j^2 = 1 \ $ then

\vspace{3mm}

\begin{equation} \label{norned sum}
\sup_{m = 1,2,\ldots} ||\sum_{j=1}^m  \gamma_j \ \eta_j \ ||B\overline{\phi} \le \theta,
\end{equation}

\vspace{3mm}

with correspondent {\it uniform} exponential tail estimate. \par

\vspace{3mm}

 \ Not analogous building may be establish this estimate   yet for the Grand Lebesgue Spaces. We will apply  the famous
 moment estimations for sums of centered independent random variables  belonging to Rosenthal, Ibragimov and Sharachmedov etc., see
\cite{Rosenthal},   \cite{Ibragimov1}, \cite{Ibragimov2}, \cite{Naimark} at all.  Namely, let $ \ \{\xi_i\}, \ i = 1,2,\ldots n \ $ be
a sequence of independent {\it centered} random variables, belonging to the space $ \ L_p(\Omega), \ p \in [2,\infty). \ $ Then

\vspace{3mm}

\begin{equation} \label{Rosenthal}
||\sum_{i=1}^n \xi_i||_p \le R(p) \ \max \left[ \  ||\sum_{i=1}^n \xi_i||_2, \ \left( \ \sum_{i=1}^n ||\xi_i||_p^p \ \right)^{1/p}  \ \right].
\end{equation}

\vspace{3mm}

 \ It is known that the optimal (minimal) value of Rosenthal "constant"  $ \ R(p)\ $ has a following form

$$
R(p)  = r \cdot \frac{p}{e \ \ln p}, \hspace{3mm} r \approx 1.77638...,
$$
see  \cite{Ibragimov1}, \cite{Ibragimov2}, \cite{Naimark}.  A slight simplification: if in addition the r.v. $ \ \{\xi_i \} \ $ are identical distributed, then

\vspace{3mm}

\begin{equation} \label{simplif}
||\sum_{i=1}^n \xi_i||_p \le R(p) \ \sqrt{n} \ ||\xi_1||_p, \ p \ge 2.
\end{equation}

\vspace{3mm}

 \ Let us consider the case when  the random variables  $ \ \xi_i  \ $ has a following form:  $ \ \xi_i =  q_i \eta_i,  \ $   where is given a vector
 $  \ \vec{q} =  \{ q_i  \}, \ i = 1,2,\ldots,n \ $ are deterministic constants, $ \  \{\eta_i \}, \ \eta := \eta_1 \ $ are centered identical distributed r.v.
 from the space $ \ L_p(\Omega). \ $ Then

\vspace{3mm}

\begin{equation} \label{Rosenth}
||\sum_{i=1}^n  \ q_i \ \eta_i||_p \le R(p) \ ||\vec{q}||_2 \ ||\eta||_p, \ p \ge 2.
\end{equation}

 \vspace{3mm}

 \ Here as ordinary $ \ ||\vec{q}||_2\stackrel{def}{=} \sqrt{ \sum_{i=1}^n q_i^2}. \ $ \par

\vspace{4mm}

\ Let again $ \ \phi(\cdot) \in \Phi \ $ be a Young - Orlicz function. Define a new following such a function (transform)

\vspace{3mm}

\begin{equation} \label{bemol token}
\phi^{(s)}(\lambda)  = \phi^{(s)}[\phi](\lambda) \stackrel{def}{=}  \sup_{v \in R} (|\lambda| |v| - \phi^*(\sqrt{|v|})).
\end{equation}

\vspace{3mm}

 \ This transform obeys a following sense. Assume that the random variable belongs to the certain  $ \ B(\phi) \ $ space and assume
 that the r.v. $ \ \kappa := \xi^2 -  \beta^2 \ = \xi^2 - \Var(\xi) \ $ satisfies also the
 famous Kramer's condition: $ \  \exists K = \const > 0 \ \rightarrow {\bf P} (|\kappa| > t) \le \exp(-Kt), \ t > 0.  \ $   Then

 \vspace{3mm}

 \begin{equation} \label{bemol sense}
 {\bf E} \exp \left[\lambda \left(\xi^2 - \beta^2 \right) \ \right] \le  e^{\phi^{(s)}(\lambda)  - \lambda \beta^2}.
 \end{equation}

\vspace{3mm}

 \ Here $ \ \lambda \in \Dom \{ \phi^{(s)} \}; \ $ herewith the case  when $ \ \lambda \in R \ $ can not be excluded.
 Note that at last $ \  (-K,K) \in \Dom \{ \phi^{(s)} \}. \ $ \par

 \ On the other words, if we introduce a new generating  $ \ \Phi \ $ function

$$
\chi(\lambda) = \chi[\phi](\lambda) \stackrel{def}{=} \chi[\phi,\beta](\lambda) := \phi^{(s)}(\lambda)  - \lambda \beta^2, \ \lambda \in \Dom \{ \phi^{(s)} \}.
$$

then
\vspace{3mm}

\begin{equation} \label{chi fun}
|| \xi^2 - \beta^2 ||B\chi \le 1.
\end{equation}

 \vspace{4mm}

   \ {\bf Remark  4.1.} \ There are close relation between tail behavior of the r.v.,  and its norms relative Grand Lebesgue Spaces as well as its norm
  relative $ \ B(\phi) \ $ ones, of course under appropriate conditions, see e.g. \cite{Koz Os Sir 4}, \cite{Os Sir 5}, \cite{Os Sir 6}, \cite{Ostrovsky7}. \par
 \ For example, let the non - zero  r.v. $ \ \xi \ $ belongs to certain $ \ B(\phi) \ $ space: $  \ ||\xi||B(\phi) \in (0,\infty). \ $ Define functions

 $$
 \beta_{\phi}(e) := \phi(\exp (y)), \ \psi_{\phi}(p) := p \exp \left(- \beta_{\phi}^*(p)/p \right).
 $$
Then
$$
 \xi  \in G\psi_{\phi}; \hspace{3mm}   ||\xi|| G\psi_{\phi} \le e^{-1} ||\xi||B \phi.
$$

\vspace{3mm}
 \ Inversely, assume the  {\it centered}  non - zero r.v. $ \ \xi \ $ belongs to some Grand Lebesgue Space $ \ G \psi: \ ||\xi||G\psi < \infty. \ $
 We suppose here that the generating function $ \ \psi(\cdot) \ $ has a following form

 $$
 \psi(p) = \psi_{\Delta}(p) = p \cdot \exp (- \Delta(p)/p ), \ p \ge 1,
 $$
where $ \ \Delta(\cdot) \ $ is  non - negative continuous convex function. Define also

$$
\phi_{\Delta}(\lambda) \stackrel{def}{=} \Delta^*(|\ln|\lambda||), \ \lambda \in R.
$$

 \ There holds

$$
||\xi||B\phi_{\Delta} \le C ||\xi|| G\psi_{\Delta}, \ C = C(\Delta) < \infty,
$$

 \vspace{3mm}

 \ {\bf Example 4.1.} For the centered r.v. $ \ \xi \ $ the subgaussian norm $ \ ||\xi||\Sub \ $ is completely equivalent to the
 following $ \ G\psi \ $ norm

$$
||\xi||G\psi_2 \stackrel{def}{=} \sup_{p \in [1,\infty)}  \left\{ \ \frac{||\xi||_p}{\sqrt{p}} \ \right\}.
$$

\vspace{3mm}

 \ {\bf Example 4.2.} For the (numerical valued) r.v. $ \ \xi \ $ the  finiteness of the following $ \ G\psi_{(\kappa)} \ $ norm

$$
||\xi||G\psi_{(\kappa)} \stackrel{def}{=} \sup_{p \ge 1} \ \left\{ \ \frac{||\xi||_p}{p^{\kappa}} \ \right\}, \ \kappa = \const \in (0,\infty),
$$
is completely equivalent to the following tail estimate

$$
T_{\xi}(c t) \le \exp \left( \ - t^{1/\kappa} \ \right), \ t \ge 1, \ \exists c = c(\kappa) \in (0,\infty).
$$

\vspace{3mm}

 \ More generally, assume that  for some r.v. $ \ \xi \ $ there holds the following estimate:  for some constants  $ \ m, c_1 > 0, \ r \in R \  $

$$
\forall p \ge 2 \ \Rightarrow \  ||\xi||_p \le c_1 \ p^{1/m} \ \ln^r p.
$$
 \ This relation is completely equivalent to the following exponential estimate for all the values $ \ t \ge e \  $

$$
T_{\xi}(t) \le \exp \left\{  \  -c_2(m,r) t^m \ \ln^{-mr}(t) \ \right\}, \exists c_2(m,r) \in (0,\infty).
$$

 \ Analogously, the relation

$$
\forall p \ge 1 \ \Rightarrow \  ||\xi||_p \le \exp \left(c_3 \ p^{\beta} \ \right),  \ \exists c_3,\beta \in (0,\infty)
$$
is equivalent to the following tail estimate:  for all the values $ \ t \ge 1 \  $

$$
T_{\xi}(t) \le \exp \left\{ \ -c_4(c_3,\beta)  \left[ \ \ln(1 + t) \  \right]^{1 + 1/\beta} \ \right\}, \  \exists c_4 = c_4(c_3,\beta) \in (0,\infty),
$$
 see e.g. \cite{Os Sir 6}, \cite{Ostrovsky7}. \par

\vspace{4mm}

 \ {\bf Remark 4.2.} \ Let us investigate the case when two random variables $ \ \xi, \eta \ $ belongs correspondingly two {\it different},
 in general case, $ \ B(\phi) \ $ spaces, say $ \ \xi \in B(\phi), \ \eta \in B(\nu), \ $   where
 $ \ \phi(\cdot), \ \nu(\cdot) \in \Phi.\ $  Define the new GLS $ \ \Phi \ $  generating  function $ \ \zeta (\cdot) \ $ as follows:

$$
\zeta (\lambda)  = \zeta[\phi,\nu](\lambda) \stackrel{def}{=} \max(\phi(\lambda), \nu(\lambda)), \ \lambda \in \Dom[\phi] \cap \Dom[\nu].
$$

 \ Then $ \ \zeta (\cdot) \in B(\phi) \ $ and herewith

$$
||\xi + \eta||B\zeta \le ||\xi||B\phi + ||\eta||B\nu.
$$

 \ Evidently, the last estimate is non - improvable. \par
  \ Quite analogously, when $ \ \xi \in G\psi, \ \eta \in G\tau, \ $ where of course  $ \ \psi(\cdot), \ \tau(\cdot) \in \Psi, \ $
   and $ \  \theta(p) := \max(\psi(p), \ \tau(p)) \ $  with $ \ p \in \Dom[\psi] \cap \Dom[\tau], \ $ then  $ \ \theta(p)  \in G\Psi \ $ and herewith

$$
||\xi + \eta||G\theta \le ||\xi||G\psi + ||\eta||G\tau.
$$

\vspace{4mm}

\section{Non - asymptotical approach.}

\vspace{3mm}

 \ Let us  return to the announced above problems.
 As we know, the r.v. $ \ \{\delta_k\} \ $ are centered asymptotically independent
 non - zero  normal identical distributed random variables. But we need
a more exact information about these variables, indeed we need estimate its Grand Lebesgue Spaces norms, in order to deduce its exponential
decreasing tails of distribution. \par

\vspace{3mm}

 \ In all these statements holds true the following expression for the empirical  Fourier coefficients

 \vspace{3mm}

\begin{equation} \label{key expression}
\hat{c}_k(n) = c_k + \frac{\delta_k}{\sqrt{n}},
\end{equation}

 \vspace{3mm}

where $ \ \{\delta_k\} \ $ are centered asymptotically independent non - zero  normal identical distributed random variables. \par

\vspace{3mm}

 \ Problem {\bf A}.  \ We suppose here that the (independent centered)  error variables $ \  \{\xi_i \}, \ i = 1,2,\ldots,n \ $
 belongs uniformly to an unit ball of the   certain $ \ B(\phi), \ \phi \in \Phi \ $ space:

$$
\exists \phi \in \Phi \ \Rightarrow \ \sup_i ||\xi_i||B(\phi)  =  1.
$$
 \ Then we deduce

$$
\sup_k \ ||\delta_k||B\overline{\phi} \le  1.
$$

  \ Following

$$
||\Sigma_1||B(\overline{\phi}) \le \frac{2}{\sqrt{n}}
$$
and

$$
||\Sigma_2||B\chi[\overline{\phi}] \le \frac{\sqrt{N}}{n}.
$$

 \ Introduce a new $ \ \Phi \ $ function

$$
\zeta(\lambda) = \zeta[\phi](\lambda) := \max \left\{ \chi[\overline{\phi}](\lambda), \overline{\phi}(\lambda)  \  \right\};
$$

Therefore

$$
|| \ \hat{\rho}_n(N) -  \rho(N)||B\zeta \le  \Theta(n,N,\rho(N)),
$$

or equally

$$
|| \ \left[ \ \hat{\rho}_n(N) -  \rho(N) \right] / \Theta(n,N,\rho(N)) \ ||B\zeta \le 1,
$$
with correspondent tail estimate. \par

\vspace{4mm}

 \ Problem {\bf B}. It is easily to calculate

$$
||\delta_k||\Sub \le \sqrt{2},
$$

as long as $ \ \max_{x \in [0,1]} |\phi_k(x)|  = \sqrt{2}. \ $ \par

\vspace{3mm}

 \ Denote

 $$
 \Sigma_1 = \frac{2}{\sqrt{n}} \cdot \sum_{k = N+1}^{2N} c_k \delta_k; \ \Sigma_2 = \frac{1}{n} \sum_{k=N+1}^{2N} (\delta_k^2 - 1);
 $$
then

$$
||\Sigma_1||_{Sub} \le  \frac{2}{\sqrt{n}} \cdot \rho(N); \ ||\Sigma_2||B \nu  \le \frac{\sqrt{N}}{n},
$$

$$
||\Sigma_1 + \Sigma_2||B \nu \le \sqrt{2} \cdot \Theta(n,N,\rho(N)).
$$

Therefore

$$
|| \ \hat{\rho}_n(N) -  \rho(N)||G\nu \le \sqrt{2} \cdot \Theta(n,N,\rho(N)),
$$

or equally

$$
|| \ \left[ \ \hat{\rho}_n(N) -  \rho(N) \right] /\left[ \ \sqrt{2} \cdot \Theta(n,N,\rho(N)) \ \right]  ||G\nu \le 1,
$$

and we deduce a non - asymptotical confidence region for $ \ \rho(N), \ 1 \le N \le n-1 \ $  by means of tail estimate
for the correspondent $ \ B\nu \ $ space.  In detail, denote for the sake of simplicity

$$
\mu = \mu (\hat{\rho}_n(N), \  \rho(N)) := \left[ \ \hat{\rho}_n(N) -  \rho(N) \right] /\left[ \ \sqrt{2} \cdot \Theta(n,N,\rho(N)) \ \right],
$$
then in accordance with (\ref{tail Bphi})

\vspace{3mm}

\begin{equation} \label{conf B}
{\bf P} \left(|\mu (\hat{\rho}_n(N), \  \rho(N))| > t \ \right)  \le 2 \exp(-\nu^*(t)), \ t > 0.
\end{equation}

\vspace{3mm}

  \ Note in addition that

 $$
 2 \exp(-\nu^*(t)) \le  \sqrt{t} \ e^{-t/2}, \ t \ge 1.
 $$

\vspace{4mm}

 \ Problem {\bf C}. Define here $ \ \Delta^2 := \Var (\xi_1); \ $  then we have in this case

 $$
 ||\delta_k|| B\upsilon \le \Delta^2.
 $$

\vspace{3mm}

 \ Therefore

\vspace{3mm}

\begin{equation} \label{upsilon}
||\Sigma_1||B \overline{\upsilon} \le \frac{ \sqrt{2} \ \Delta^2}{\sqrt{n}} \cdot \rho(N).
\end{equation}

\vspace{3mm}

 \ Further, here $ \ \Sigma_2  = n^{-1} \sum_{k=N+1}^{2N} (\delta_k^2 - \Delta^2),  \ $ and we find

 \vspace{3mm}

\begin{equation} \label{sigma2}
||\Sigma_2||B \left( \overline{\chi[\phi_2]}  \ \right)  \le \sqrt{2} \cdot \frac{\sqrt{N}}{n} \cdot \Delta^2.
\end{equation}

\vspace{3mm}

 \ We conclude as above

$$
||  \left\{ \ \left[ \ \hat{\rho}_n(N) -  \rho(N) \right] /\left[ \ \sqrt{2} \cdot \Theta(n,N,\rho(N)) \ \right]\ \right\} || B \left( \overline{\chi[\phi_2]} \ \right) \le  \Delta^2,
$$
  with accordance  correspondent tail estimation (\ref{tail Bphi}).\\

\vspace{3mm}
 \hspace{3mm} {\bf Examples 5.1.} Imposed conditions on the function $ \ \rho(N) \ $ are satisfied for instance when
 $ \  \rho(N) \sim  \ N^{-\alpha} L(N), \  \alpha = \const >  0,  \ $ i.e. regular varying function,
  where $ \ L = L(N), \ N = 1,2,\ldots \ $ is positive finite
 slowly varying function at infinity,  or when $ \ \rho(N) \asymp q^N, \ q = \const \in (0,1). \ $ \par
 \ For instance,

\vspace{3mm}

\begin{equation} \label{alpha}
\rho(N) = c_1 \ N^{-\alpha} \ \ln^{\gamma}(N +1),  \ c_1, \ \alpha  = \const >  0, \ \gamma = \const \ge 0,
\end{equation}

 \vspace{3mm}
{\it quasi - power model,}  or may be
 \vspace{3mm}
\begin{equation} \label{qexpression}
\rho(N) = c_2 \ N^{\kappa} \ q^N, \ c_2 = \const \in (0,\infty), \kappa = \const, \ q = \const \in (0,1),
\end{equation}
{\it quasi exponential model.}\par
 \vspace{3mm}

 \ {\bf Remark 5.1.} \ Notable that was made many success  computer simulations  of all three considered non - parametrical statistical
 problem for these examples of the key functions $ \ \rho(N). \ $ \par

\vspace{3mm}

 \ {\bf Remark 5.2. Functional approach.} \ Let us consider the following confidence probability
 \vspace{3mm}
$$
Q = Q_{a,b} (\vec{v}) \stackrel{def}{=} {\bf P} \left(\ \cup_{N \in [a,b]} \left\{ \ |\mu(\hat{\rho}_n(N),\rho(N))| > v(N) \ \right\} \right).
$$
\vspace{3mm}
 \ Here $ \  [a,b], \ 1 \le a < b \le n-1 \ $ is an integer subinterval of $  \ [1, n-1] \ $ one, $ \ \vec{v} = \{ v(a), v(a+1),  \ldots, v(b) \} \ $
is arbitrary positive numerical vector.  We deduce obviously

$$
Q \le \sum_{N = a}^b {\bf P} \left(\ |\mu(\hat{\rho}_n(N),\rho(N))| > v(N) \ \right).
$$
 \vspace{3mm}

  \ We have in the particular case when $ \ v(N) = w = \const > 0  \ $ the following maximal distribution estimation for the confidence region
  in the  so - called  "uniform norm"

$$
{\bf P} \left( \  \max_{N \in [a,b]}  |\mu(\hat{\rho}_n(N),\rho(N))| > w  \ \right) \le \sum_{N = a}^b {\bf P} \left(\ |\mu(\hat{\rho}_n(N),\rho(N))| > w \ \right).
$$

\vspace{3mm}

 \ {\bf Remark 5.3.} \ The offered  before  {\it the trajectory} estimate $ \ \hat{\rho}_n(N), \ 1 \le N \le n-1 \ $ may be used, for example, for the consistent
 its  {\it parametric} approximation, see e.g. (\ref{alpha}):

 \vspace{3mm}

 \begin{equation} \label{approximation}
 \hat{\rho}_n(N) \approx c_1 \ N^{-\alpha} \ \ln^{\gamma}(N +1),  \ c_1, \ \alpha  = \const >  0, \ \gamma = \const \ge 0,
 \end{equation}

\vspace{3mm}
by means of the famous Minimum Square Estimation  (MSE) method. \par
 \ In detail, introduce the following function

 \vspace{3mm}

\begin{equation} \label{Subject opt pol}
Z = Z_n(c_1, \alpha, \gamma) \stackrel{def}{=}  \sum_{N=1}^{n-1} \left[ \ \hat{\rho}_n(N) - c_1 N^{-\alpha} \ln^{\gamma}(N+1) \  \right]^2,
 \end{equation}

 \vspace{3mm}

where as before  $ \ c_1, \ \alpha >  0, \ \gamma \ge 0. \ $ The consistent as $ \ n \to \infty \ $ estimates  $ \ (\hat{c_1},  \ \hat{\alpha},  \hat{\gamma}) \ $
of the parameters $ \ (c_1, \alpha, \gamma) \ $ are the point os minimum of introduced function

\vspace{3mm}

\begin{equation} \label{point min}
(\hat{c_1},  \ \hat{\alpha},  \hat{\gamma}) = \argmin_{ \{c_1,  \ \alpha, \ \gamma \} }  Z_n(c_1, \alpha, \gamma).
\end{equation}

\vspace{3mm}

 -  the quasi - power model.     Analogously for the  quasi - exponential  model we get to the following extremal problem

\vspace{3mm}

\begin{equation} \label{qexpression}
(\hat{c_2},  \hat{\kappa}, \hat{q}) := \argmin_{ c_2,\kappa, q} \sum_{N=1}^{n-1} \left[\hat{\rho}_n(N)  - c_2 \ N^{\kappa} \ q^N \ \right]^2.
\end{equation}

\vspace{4mm}

 \ {\bf Remark 5.4.} Both the problems  (\ref{point min})  and (\ref{qexpression})  lead us to the complicated optimization problems.
All the  iterative numerical methods (Newton etc.) need an initial approximation. It may be found from the simplification of offered
model, of  course, by means of logarithmical "transform," for example

\vspace{3mm}

 \begin{equation} \label{logarithization}
 \log \hat{\rho}_n(N) \approx  c_3 - \alpha_1 \cdot \ln N  \ + \gamma_1  \cdot \ln\ln(N +1).
 \end{equation}

\vspace{3mm}

 \ Let's form an auxiliary loss  function

\vspace{3mm}

$$
W = W(c_3, \alpha_1, \gamma_1) := \sum_{N=1}^{n-1} \left[ \  \ln \hat{\rho}_n(N) - c_3 + \alpha_1 \cdot \ln N - \gamma_1 \ln \ln N \ \right]^2.
$$

\vspace{3mm}

 \ The problem of minimization of this function

 \vspace{3mm}

$$
\left\{ \ \hat{c_3},  \hat{\alpha_1}, \hat{\gamma_1} \right\} := \argmin_{c_3, \alpha_1, \gamma_1} W(c_3, \alpha_1, \gamma_1)
$$
 leads to the linear system equations. The values

$$
c_1^0 = \exp \hat{c_3}, \ \alpha_1^0 = \hat{\alpha_1}, \ \gamma_1^0 = \hat{\gamma_1}
$$
may be used as an initial approximation for the problem (\ref{point min}). \par

 \vspace{3mm}

  \ The case of the quasy - exponential approximation may be considered quite analogously. \par

\vspace{4mm}

\section{Concluding remarks.}

\vspace{4mm}

 \hspace{3mm} It is interest in our opinion to generalize obtaining results into  multivariate and other
 non - parametrical statistical problems, as well as into the more general domains as a standard interval $ \ [0,1]. $ \par

\vspace{6mm}

\vspace{0.5cm} \emph{Acknowledgement.} {\footnotesize The first
author has been partially supported by the Gruppo Nazionale per
l'Analisi Matematica, la Probabilit\`a e le loro Applicazioni
(GNAMPA) of the Istituto Nazionale di Alta Matematica (INdAM) and by
Universit\`a degli Studi di Napoli Parthenope through the project
\lq\lq sostegno alla Ricerca individuale\rq\rq .\par


\vspace{6mm}

\begin{thebibliography}{79}


\bibitem{Ahmed Fiorenza Formica at all}
{\bf I. Ahmed, A. Fiorenza, M.R. Formica, A. Gogatishvili, J.M. Rakotoson.}
{\it Some new results related to Lorentz G-Gamma spaces and interpolation.}
J. Math. Anal. Appl., 483, {\bf 2}, (2020).


\bibitem{anatriellofiojmaa2015}
{\bf G.~Anatriello} and {\bf A.~Fiorenza.} {\it Fully measurable
grand Lebesgue spaces}. J. Math. Anal. Appl. \textbf{422} (2015),
no.~2, 783--797.

\bibitem{anatrielloformicaricmat2016}
{\bf G.~Anatriello} and {\bf M.~R.~Formica.} {\it Weighted fully
measurable grand Lebesgue spaces and the maximal theorem}. Ric. Mat.
\textbf{65} (2016), no.~1, 221--233.

\bibitem{Bobrov}
{\bf Bobrov P.B., Ostrovsky E.I.} {\it  Adaptive estimations of regression, density and spectrums.}  1996, Laboratory, V.5, pp. 55 -- 58.

\bibitem{Bobrov2}
{\bf Bobrov P.B., Ostrovsky E.I.} {\it  Adaptive non - parametrical optimal  estimations of regression, density and spectrums.}  1997,
Problems of information transmission, V.4, 33 - 41.

\bibitem{Buld Koz AMS}
{\bf Buldygin V.V., Kozachenko Yu.V. }  {\it Metric Characterization of Random
 Variables and Random Processes.} 1998, \ Translations of Mathematics Monograph, AMS, v.188.



\bibitem{Buldygin-Mushtary-Ostrovsky-Pushalsky} {\bf V.~V.~Buldygin, D.~I.~Mushtary, E.~I.
~Ostrovsky} and {\bf M.~I.~Pushalsky.} {\it New Trends in
Probability Theory and Statistics.} Mokslas (1992), V.1, 78--92;
Amsterdam, Utrecht, New York, Tokyo.

\bibitem{caponeformicagiovanonlanal2013}
{\bf C.~Capone, M.~R.~Formica} and {\bf R.~Giova.} {\it Grand
{L}ebesgue spaces with respect to measurable functions}. Nonlinear
Anal. \textbf{85} (2013), 125--131.



\bibitem{Ermakov etc. 1986}
{\bf S. V. Ermakov, and E. I. Ostrovsky.} {\it Continuity Conditions, Exponential Estimates, and the Central Limit Theorem for Random Fields.}
 Moscow, VINITY,  1986. (in Russian).



\bibitem{Fiorenza1} {\bf A.~Fiorenza.} {\it Duality and reflexivity in grand Lebesgue
spaces.} Collect. Math. \textbf{51} (2000), no. 2, 131--148.

\bibitem{Fiorenza4}
{\bf A.~Fiorenza} and {\bf G.~E.~Karadzhov.} {\it Grand and small
Lebesgue spaces and their analogs}, Z. Anal. Anwendungen \textbf{23}
(2004), no.~4, 657--681.

\bibitem{fioguptajainstudiamath2008}
{\bf A.~Fiorenza, B.~Gupta} and {\bf P.~Jain.} {\it The maximal
theorem for weighted grand Lebesgue spaces}. Studia Math.
\textbf{188} (2008), no.~2, 123--133.


\bibitem{Fiorenza-Formica-Gogatishvili-DEA2018}
{\bf A.~Fiorenza, M.~R.~Formica} and {\bf A.~Gogatishvili.} {\it On
grand and small Lebesgue and Sobolev spaces and some applications to
PDE's}. \emph{Differ. Equ. Appl.} \textbf{10} (2018), no.~1, 21--46.

\bibitem{fioforgogakoparakoNAtoappear}
{\bf A.~Fiorenza, M. R.~Formica, A.~Gogatishvili, T.~Kopaliani} and
{\bf J.~M. Rakotoson.} {\it Characterization of interpolation
between grand, small or classical Lebesgue spaces}. Nonlinear Analysis, Vol.177, Part {\bf B,} Dezember 2018,
pages 422 \ - \ 453.


\bibitem{fioformicarakodie2017}
{\bf A.~Fiorenza, M.~R.~Formica} and {\bf J.~M. Rakotoson.} {\it
Pointwise estimates for {$ \ G\Gamma$c \ }-functions and applications}.
Differential Integral Equations, \textbf{30}, (2017), no. \ ~11 \ - \ 12,
809 \ - \ 824.

\bibitem{formicagiovamjom2015}
{\bf M.~R. Formica} and {\bf R.~Giova.} {\it Boyd indices in
generalized grand Lebesgue spaces and applications}. Mediterr. J.
Math., \textbf{12}, (2015), no.~3, 987 \ - \ 995.



\bibitem{Golubev Lepsky Levit}
{\bf G. Golubev, O.Lepsky. B.Levit} {\it On Adaptive Estimation Using the Sup-Norm Losses.}\\
  Mathematical Methods of Statistics, May 2000,  V. 2; pp.34 - 51.
.



\bibitem{Ibragimov1}
{\bf Ibragimov R., Sharachmedov Sh.} {\it On the exact constant in the Rosenthal
Inequality.} Theory Probab. Appl., 1997, V. {\bf 42,}  p. 294 - 302.


\bibitem{Ibragimov2}
{\bf Ibragimov R., Sharachmedov Sh.}  {\it The exact constant in the Rosenthal Inequality for sums random variables
 with mean zero.} Probab. Theory Appl., 2001, V. {\bf 46,} \ 1, \ p. 127 - 132.

\bibitem{Iwaniec2}
{\bf T.~Iwaniec} and {\bf C.~Sbordone.} {\it On the integrability of
the Jacobian under minimal hypotheses.} Arch. Rational Mech. Anal.
\textbf{119} (1992), no.~2, 129--143.

\bibitem{KozOs}
{\bf Yu.~V.~Kozachenko} and {\bf E.~I.~Ostrovsky.} {\it The Banach
Spaces of random variables of sub-Gaussian type.} of Probab. and
Math. Stat., \textbf{32} (1985), (in Russian). Kiev, KSU, 43--57.


\bibitem{Kozachenko at all 2018}
{\bf Yu.V. Kozachenko, Yu.Yu. Mlavets, and N.V. Yurchenko.} {\it Weak convergence of stochastic processes from spaces} $F_\psi(\Omega).$
STATISTICS, OPTIMIZATION AND INFORMATION COMPUTING, Vol.6,  June 2018, pp. 266 - 277.

\bibitem{Koz Os Sir 4}
{\bf Kozachenko Yu.V., Ostrovsky E., Sirota L.} {\it Relations between exponential tails, moments and
moment generating functions for random variables and vectors.} \\
 arXiv:1701.01901v1 [math.FA] 8 Jan 2017


\bibitem{liflyandostrovskysirotaturkish2010}
{\bf E.~Liflyand, E.~Ostrovsky} and {\bf L.~Sirota.} {\it Structural
properties of bilateral grand {L}ebesque spaces}. Turkish J. Math.
\textbf{34} (2010), no.~2, 207--219.


\bibitem{Naimark}
{\bf Naimark B., Ostrovsky E.} {\it Exact Constants in the Rosenthal Moment Inequalities for Sums of
independent centered Random Variables.} \\
arXiv:math/0411614 [math.PR]


\bibitem{Ostrovsky1}
{\bf E.Ostrovsky.} {\it Exponential estimates for random fields and
its applications.} 1999, OINPE, Moscow - Obninsk.


 \bibitem{Ostrovsky2}
 {\bf E.Ostrovsky.} {\it Exponential estimate in the Law of Iterated
 Logarithm in Banach Space.} Math. Notes \textbf{56} (1994), no.~ 5-6, 1165--1171.

\bibitem{Ostrovsky3}
 {\bf Ostrovsky E. and Sirota L.} {\it Moment Banach spaces: theory and applications.}
HIAT Journal of Science and Engineering, C, Volume 4, Issues 1 - 2, pp. 233 \ - \ 262, (2007).


 \bibitem{Ostrovsky4}
 {\bf  Ostrovsky E., Sirota L.} {\it Prokhorov-Skorokhod continuity of random fields.
 A natural approach.} \\ arXiv:1710.05382v1 [math.PR] 15 Oct 2017

\bibitem{Os Sir 5}
{\bf Ostrovsky E., Sirota L.}  {\it Moment Banach Spaces: Theory and Applications.} \\
HIAT Journal of Science and Engineering, Holon, Israel, v. 4, Issue 1-2, (2007), 233 \ - \ 262.


\bibitem{Os Sir 6}
{\bf Ostrovsky E., Sirota L.}  {\it Sharp moment and exponential tail estimates for U - statistics.}
arXiv:1602.00175 [math.ST].


\bibitem{Ostrovsky7}
{\bf Ostrovsky E.} {\it Bide-side exponential and moment inequalities for tails of
distribution of Polynomial Martingales.} \\
arXiv: math.PR/0406532 v.1 Jun. 2004

\bibitem{Rosenthal}
{\bf Rosenthal H.P.} {\it On the subspaces of} $ \ Lp, (p > 2) \ $ {\it spanned by sequences of
independent variables. } Israel J. Math., 1970, N o3, p. 273 - 303.

\bibitem{Samko-Umarkhadzhiev}
{\bf S.~G.~Samko} and {\bf S.~M.~Umarkhadzhiev.} {\it On
Iwaniec-Sbordone spaces on sets which may have infinite measure.}
Azerb. J. Math. \textbf{1} (1)\ (2011), 67--84.

\bibitem{Samko-Umarkhadzhiev-addendum}
{\bf S.~G.~Samko} and {\bf S.~M.~Umarkhadzhiev.} {\it On
Iwaniec-Sbordone spaces on sets which may have infinite measure:
addendum.} Azerb. J. Math \textbf{1} (2) \ (2011), 143--144.

\bibitem{Tchentzov1}
{\bf Tchentzov N.N.} {\it Density estimation based on the sample.} Doklady of the
Soviet Akademy of Science. 1962,  V. 147, Issue  1, pp.  45 - 48. (In Russian).


\bibitem{Tchentzov2}
{\bf Tchentzov N.N.} {\it Statistical decision rules and optimal conclusions.} Moscow,
Nauka, 1972, 520 pp., (In Russian).

\bibitem{Yudovich1}
{\bf V. I. Yudovich.} {\it Non-stationary flows of an ideal incompressible fluid.} Z. Vychisl. Mat.
i Mat. Fiz., 3: 1032 \ - \ 1066, (Russian), 1963.

\bibitem{Yudovich2}
{\bf V. I. Yudovich.} {\it Uniqueness theorem for the basic nonstationary problem in the dynamics of an ideal
incompressible fluid.} Math. Res. Lett., 2(1):27�38, 1995. 1, 3, 4,5.

\end{thebibliography}
\end{document}